      \documentstyle[12pt, righttag]
{amsart}

    \input amssym.def
    \input amssym

    \newtheorem{propo}{Proposition}
    \newtheorem{theo}[propo]{Theorem}
    
    \newtheorem{corol}[propo]{Corollary}

    \theoremstyle{definition}
    \newtheorem{defi}[propo]{Definition}

    \theoremstyle{remark}
    \newtheorem{rema}[propo]{Remark}

\numberwithin{equation}{section}
\numberwithin{propo}{section}

\begin{document}
    \title[Intertwining operator algebras and
conformal field theories]
{Intertwining Operator Algebras, Genus-Zero Modular Functors and
Genus-Zero Conformal Field
Theories}
    \author{Yi-Zhi Huang}
    \address{Department of Mathematics, Rutgers University,
New Brunswick, NJ 08903}
    \email{yzhuang@@math.rutgers.edu}
    \thanks{Supported in part by NSF grants
DMS-9301020  and DMS 9596101 and by DIMACS, an
NSF Science and Technology Center funded under contract STC-88-09648.}
	\bibliographystyle{alpha}
	\maketitle

	\newcommand{\nno}{\nonumber}
	\newcommand{\lbar}{\bigg\vert}
\newcommand{\mbar}{\mbox{\huge $\vert$}}
	\newcommand{\p}{\partial}
	\newcommand{\dps}{\displaystyle}
	\newcommand{\bra}{\langle}
	\newcommand{\ket}{\rangle}
 \newcommand{\res}{\mbox{\normalshape Res}}
 \newcommand{\epf}{\hspace{2em}$\Box$}
 \newcommand{\epfv}{\hspace{1em}$\Box$\vspace{1em}}
\newcommand{\nord}{\mbox{\scriptsize ${\circ\atop\circ}$}}
\newcommand{\wt}{\mbox{\normalshape wt}\ }
\newcommand{\clr}{\mbox{\normalshape clr}\ }

\begin{abstract}
We describe the construction of the genus-zero parts
of conformal field theories in the sense of G. Segal
{}from representations of
vertex operator algebras satisfying certain conditions.
The construction is divided into four steps and each step gives a
mathematical structure of independent interest. These mathematical
structures are intertwining operator algebras, genus-zero modular
functors, genus-zero holomorphic weakly conformal field theories, and
genus-zero  conformal field theories.
\end{abstract}

\section{Introduction}

Conformal field theory (see for example, \cite{BPZ}, \cite{W} and
\cite{MS})
 is a physical theory related to many branches of
mathematics, e.g., infinite-dimensional Lie algebras and Lie groups,
sporadic finite simple groups, modular forms and modular functions,
elliptic genera and elliptic cohomology, Calabi-Yau manifolds
and mirror symmetry, and quantum groups and $3$-manifold invariants.
Recently there are efforts by mathematicians to develop
conformal field theory as a rigorous mathematical theory so that it
can be used seriously in the near future
to solve mathematical problems. A mathematical
precise
definition of conformal field theory was first given by Segal \cite{S}
several years ago.
Segal's definition is geometric and
is motivated by the path integral formulation of some
conformal field theories. Though this definition is conceptually simple,
it is difficult to see {}from it
some very subtle but important properties of conformal
field theories. Because of this,
it is even more difficult to construct directly conformal
field theories satisfying Segal's geometric definition.

In practice, physicists and mathematicians  have studied concrete models of
conformal field theories for many years and the methods used are mostly
algebraic. These studies can also be summarized to give algebraic
formulations of conformal field theory. The theory of vertex operator
algebras is such an algebraic formulation.
Many
concrete examples of conformal field theories are formulated and
studied in algebraic formulations and
very detailed calculations can be carried out in these formulations.
But on the other hand, algebraic formulations have a disadvantage that
the higher-genus parts of conformal field theories cannot be even formulated,
and thus in algebraic formulations,
it is not easy to see the geometric and topological applications of
conformal field theory .

Intuitively, it is expected that  algebraic and geometric formulations
be equivalent. But  this equivalence, especially, the construction
of conformal field theories satisfying the geometric definition {}from
conformal field theories satisfying the algebraic definition is a highly
nontrivial  mathematical problem.
In the present paper, we describe the construction of the genus-zero parts
of conformal field theories
{}from representations of
vertex operator algebras satisfying certain conditions.

This construction can actually be divided into four steps and each
step gives us a mathematical structure of independent interest.  The
first step is to construct an ``intertwining operator algebra'' {}from
the irreducible representations of a vertex operator algebra satisfying certain
 conditions.  The second step is to construct a ``genus-zero
modular functor'' (a partial operad of a certain type) {}from an
intertwining operator algebra. The third step is to construct a
``genus-zero holomorphic weakly conformal field theory'' (an algebra
over the partial operad of the genus-zero modular functor satisfying
certain additional properties) {}from the intertwining operator
algebra. The last step is to construct a ``genus-zero conformal field
theory'' {}from a genus-zero holomorphic weakly conformal field theory
when the genus-zero holomorphic weakly conformal field theory is
unitary in a certain sense.  The first and the second steps are
described in Sections 3 and 4, respectively. The third and the fourth
are both described in Section 5.

We introduce the notions of intertwining operator algebra and
genus-zero modular functor in Section 3 and 4, respectively. The
notions of genus-zero holomorphic weakly conformal field theory and
genus-zero conformal field theory are introduced in Section 5.  The
notion of intertwining operator algebra is new. The other three
notions are modifications of the corresponding notions introduced by
Segal \cite{S} in the genus-zero case. We also give a brief description
of the operadic formulation of the notion of vertex operator algebra in
Section 2.

The construction described in this paper depends on the solutions of
two problems. The first is the precise geometric description of
the central charge of a vertex operator algebra. This problem is solved
completely in \cite{H4} (see also \cite{H1} and \cite{H2}).
The second is the associativity
of intertwining operators for a vertex operator algebra satisfying
certain conditions. This is proved in \cite{H3} using the
tensor product theory for modules for a vertex operator algebra
(see \cite{HL3}--\cite{HL6}, \cite{H3}). All the results described in
the present paper are consequences of the results obtained in \cite{H1},
\cite{H4}, \cite{HL3}, \cite{HL4}, \cite{HL6}, \cite{H3} and \cite{H5}.

\subsection*{Acknowledgment}
I would like to thank J.-L. Loday,
J. Stasheff and A. A. Voronov for inviting me to attend these
two conferences.

\vspace{1em}
\noindent {\bf Notations}:

\noindent $\Bbb{C}$: the (structured set of) complex numbers.

\noindent $\Bbb{C}^{\times}$: the nonzero complex numbers.

\noindent $\Bbb{R}$: the real numbers.

\noindent $\Bbb{Z}$: the integers.

\noindent $\Bbb{Z}_{+}$: The positive integers.

\noindent $\Bbb{N}$: the nonnegative integers.

\section{Operadic formulation of the notion of vertex operator algebra}

We give a brief description of the operadic formulation
of the notion of vertex operator algebra in this section. We assume that the
reader is familiar with operads and algebras over operads. See, for example,
\cite{M}
for an introduction to these notions.
The material in this section
is taken {}from  \cite{H1},
\cite{HL1},
\cite{HL2} and \cite{H4}. Since these references contain detailed and precise
definitions and theorems, here we only present  the
main ideas. The main purposes of this section is to convince those readers
unfamiliar with vertex operator algebras that vertex
operator algebras are in fact very natural mathematical objects, and to
introduce informally  some basic concepts and notations in preparation for
later sections. For details, see the references above. Note that though
the operadic formulation
of the notion of vertex operator algebra described below is very
natural {}from the viewpoint of operads and geometry, historically,
vertex operator algebras occurred in mathematics and physics for completely
different reasons.
See the introduction of
\cite{FLM} for a detailed historical discussion.  For the algebraic
formulation of the notion of vertex (operator) algebra, see \cite{B},
\cite{FLM} and \cite{FHL}.

As a motivation, we begin with associative algebras.
Let $C(j)$, $j\in \Bbb{N}$, be the moduli space of
circles (i.e., compact connected smooth one-dimensional manifolds)
with $j+1$ ordered points (called {\it punctures}), the zeroth
negatively oriented, the others positively oriented, and with smooth local
coordinates vanishing at these punctures.  Then it is easy to see that
$C(j)$ can be identified naturally with the set of permutations
$(\sigma(1), \dots, \sigma(j))$ of $(1,
\dots, j)$. Since this set can in turn be identified
in the obvious way with the symmetric group $S_{j}$, the moduli space
$C(j)$ can also be naturally identified with $S_{j}$, with the group
$S_{j}$ acting on $C(j)$ according to the left multiplication action
on $S_{j}$. That is, $S_{j}$ permutes the orderings of the positively
oriented punctures. Given any two circles with punctures and local
coordinates vanishing at the punctures, we can define the notion of
sewing them together at any positively oriented puncture on the first
circle and the negatively oriented puncture on the second circle by
cutting out an open interval of length $2r$ centered at the positively
oriented puncture on the first circle (using the local coordinate) and
cutting out an open interval of length $2/r$ around the negatively
oriented puncture on the second circle in the same way, and then by
identifying the boundaries of the remaining parts using the two local
coordinates and the map $t\longmapsto 1/t$; we assume that the
corresponding closed intervals contain no other punctures.  The
ordering of the punctures of the sewn circle is given by ``inserting''
the ordering for the second circle into that for the first.  Note that
in general not every two circles with punctures and local coordinates
can be sewn together at a given positively oriented puncture on the
first circle.  But it is clear that such a pair of circles with
punctures and local coordinates is equivalent to a pair which can be
sewn together. Also, the equivalence class of the sewn circle with
punctures and local coordinates depends only on the two equivalence
classes.  Thus we obtain a sewing operation on the moduli space of
circles with punctures and local coordinates. Given an element of
$C(k)$ and an element of $C(j_{s})$ for each $j_{s}$, $s=1, \dots,k$,
we define an element of $C(j_{1}+\cdots +j_{k})$ by sewing the element
of $C(k)$ at its $s$-th positively oriented puncture with the element
of $C(j_{s})$ at its negatively oriented puncture, for $s=1,\dots, k$,
and by ``inserting'' the orderings for the elements of the $C(j_s)$
into the ordering for the element of $C(k)$.  The identity is the
unique element of $C(1)$.
It is straightforward to verify that $C=\{C(j)\}_{j\in \Bbb{N}}$ is
indeed an operad.

Since in one dimension smooth structures and conformal structures are
the same, the moduli space $C(j)$ can also be thought of as the moduli
space of circles with conformal structures, with $j+1$ ordered
punctures with the zeroth negatively oriented and the others
positively oriented, and with local conformal coordinates vanishing at
these punctures. In fact, we have just defined three operads, namely,
$C$, the corresponding conformal moduli space and $\{ S_{j}\}_{j\in
\Bbb{N}}$, and these three operads are isomorphic.

It is well known that the category of algebras over the operad $\{
S_{j}\}_{j\in \Bbb{N}}$ is isomorphic to the category of associative
algebras. Thus the category of algebras over $C$ (or the category of
algebras over the operad of the corresponding  conformal moduli space) is
isomorphic to the category of associative algebras.

We have just seen that associative algebras are algebras over an operad
obtained {}from one-dimensional objects. It is very natural to consider
algebras over operads obtained {}from two-dimensional objects. The
simplest two-dimensional objects are topological spheres, i.e.,
genus-zero compact connected smooth two-dimensional manifolds. If we
consider the operad of the moduli spaces of genus-zero compact
connected smooth two-dimensional manifolds with punctures and local
coordinates vanishing at these punctures, then the category of
algebras over this operad is isomorphic to the category of commutative
associative algebras.  We do not obtain new algebras. On the other
hand, in dimension two, conformal structures are much richer than
smooth structures. Since conformal structures in two real dimensions
are equivalent to complex structure in one complex dimension, it is
natural to consider operads constructed {}from genus-zero compact
connected one-dimensional complex manifolds with punctures and local
coordinates vanishing at these punctures.

Roughly speaking, vertex operator algebras are algebras over a
certain operad constructed {}from
genus-zero compact connected one-dimensional complex
manifolds with punctures and local
coordinates. But the operad in this case has an analytic
structure and this analytic structure becomes algebraic in a certain
sense when the operad is extended to a partial operad. Vertex
operator algebras have properties which are
 not only the reflections of the operad structure above, but also the
reflections of the analytic structure
of the extended partial operad. These analyticity
properties of vertex operator algebras are
very subtle and important features of
the theory of vertex operator algebras. They are essential to the
constructions and many applications of conformal field theories.

We now begin to describe the notion of vertex operator algebra using the
language of operads.
A {\it
sphere with $1+n$ tubes} ($n\in \Bbb{N}$) is a  genus-zero compact
connected one-dimensional complex manifold $S$ with $n+1$ distinct,
ordered points $p_{0}, \dots , p_{n}$ (called {\it punctures}) on
$S$ with $p_{0}$ negatively oriented and the other punctures
positively oriented, and with local analytic coordinates $(U_{0},
\varphi _{0}), \dots, (U_{n}, \varphi _{n})$ vanishing at the
punctures $p_{0}, \dots , p_{n}$, respectively, where for $i=0,
\dots, n$, $U_{i}$ is a local coordinate neighborhood at $p_{i}$
(i.e., an open set containing $p_{i}$) and $\varphi _{i}: U_{i}
\to
\Bbb{C}$,
satisfying $\varphi _{i}(p_{i})=0$, is a local analytic coordinate map
vanishing at $p_{i}$.  Let $S_{1}$ and $S_{2}$ be spheres with $1+m$
and $1+n$ tubes, respectively. Let $p_{0},
\dots, p_{m}$ be the punctures of $S_{1}$, $q_{0},\dots,
q_{n}$ the punctures of $S_{2}$, $(U_{i}, \varphi_{i})$ the local
coordinate at $p_{i}$ for some fixed $i$, $0<i\le m$, and $(V_{0},
\psi_{0})$ the local coordinate at $q_{0}$. Assume that there exists a
positive number $r$ such that $\varphi _{i}(U_{i})$ contains the
closed disc $\bar{B}_{0}^{r}$ centered at $0$ with radius $r$ and
$\psi_{0}(V_{0})$ contains the closed disc $\bar{B}_{0}^{1/r}$
centered at $0$ with radius $1/r$.  Assume also that $p_{i}$ and
$q_{0}$ are the only punctures in $\varphi_{i}^{-1}(\bar{B}_{0}^{r})$
and $\psi _{0}^{-1}(\bar{B}_{0}^{1/r})$, respectively. In this case we
say that {\it the $i$-th tube of $S_{1}$ can
be sewn with the zeroth tube of $S_{2}$}. {}From
$S_{1}$ and $S_{2}$, we obtain a sphere with $1+(m+n-1)$
tubes by cutting $\varphi_{i}^{-1}(B_{0}^{r})$ and
$\psi_{0}^{-1}(B_{0}^{1/r})$ {}from $S_{1}$ and $S_{2}$, respectively,
and then identifying the boundaries of the resulting surfaces using
the map $\varphi_{i}^{-1} \circ J \circ
\psi_{0}$ where $J$ is the map {}from $\Bbb{C}^{\times}$ to
itself given by $J(z)=1/z$. The negatively oriented puncture of
this sphere with tubes is $p_{0}$ and the positively
oriented punctures (with ordering) of this
sphere 	with tubes are $p_{0}, \dots,p_{i-1}$, $q_{1}, \dots,
q_{n}$, $p_{i+1}, \dots,p_{m}$. The local coordinates vanishing at
these punctures are given in the obvious way. Thus we have a partial
operation---the {\it sewing operation}---in the collection of
spheres with tubes. We define the notion of conformal equivalence between two
spheres with tubes in the obvious way except that  two
spheres with tubes are also said to be conformally equivalent if the only
differences between them are local coordinate neighborhoods
at the punctures. For any $n\in \Bbb{N}$,
the space of equivalence
classes of spheres with $1+n$ tubes is called the {\it moduli space of
spheres with $1+n$ tubes}.  For $n\in \Bbb{Z}_{+}$,
the moduli space of spheres with $1+n$ tubes
can be identified with $K(n)=M^{n-1}\times H \times
(\Bbb{C}^{\times} \times H)^{n}$ where $M^{n-1}$ is the set of elements of $
\Bbb{C}^{n-2}$ with nonzero and distinct components and
$H$ is the set of all sequences $A$ of complex
numbers such that $(\mbox{\rm exp}(
\sum_{j>0}A_{j}x^{j+1}\frac{d}{dx}))x$ is a
convergent power series in some neighborhood of $0$. We think of each
element of $K(n)$, $n\in \Bbb{Z}_{+}$, as the sphere $\Bbb{C}\cup \{ \infty\}$
equipped with negatively oriented puncture  $\infty$ and
positively oriented ordered punctures $z_{1},\dots,z_{n-2}$, $0$,
with an element of $H$ specifying the local coordinate at $\infty$ and
with $n$ elements of $\Bbb{C}^{\times} \times H$
specifying the local coordinates at the
other punctures. Analogously, the moduli space of spheres with $1+0$
tube can be identified with $K(0)=\{B\in H\;|\;B_{1}=0\}$.
{}From now on we shall refer to $K(n)$ as the moduli
space of spheres with $1+n$ tubes, $n\in \Bbb{N}$.
The sewing operation for spheres with
tubes induces a partial operation on the $\cup_{n\in \Bbb{N}}K(n)$.
It is still called the
sewing operation and is denoted $_{^{i}}\infty_{^{0}}$.

Let $K=\{K(n)\}_{n\in \Bbb{N}}$. For any $k\in \Bbb{Z}_{+}$, $j_{i},
\dots, j_{k}\in \Bbb{N}$, and any elements of $K(k)$, $K(j_{1}),
\dots, K(j_{k})$, by successively sewing the $0$-th tube of elements
of $K(j_{i})$, $i=1, \dots, k$, with the $i$-th tube of the element of
$K(k)$, we obtain an element of $K(j_{1}+\cdots +j_{k})$ if the
conditions to perform the sewing operation are satisfied. Thus we
obtain a partial map $\gamma_{K}$ {}from $K(k)\times K(j_{1})\times
\cdots K(j_{k})$ to $K(j_{1}+\cdots +j_{k})$.  Let $I\in K(1)$ be
equivalence class of the standard sphere $\Bbb{C}\cup \{ \infty\}$
with $\infty$ the negatively oriented puncture, $0$ the only
positively oriented puncture, and with standard local coordinates
vanishing at $\infty$ and $0$. For any $j\in \Bbb{N}$, $\sigma\in
S_{j}$ and $Q\in K(j)$, $\sigma(Q)$ is defined to be the conformal
equivalence class of spheres with $1+j$ tubes obtained {}from members of
the class $Q$ by permuting the orderings of their positively oriented
punctures using $\sigma$. Thus $S_{j}$ acts on $K(j)$. It is easy to
see that $K$ together with $\gamma_{K}$, $I$ and the actions of
$S_{j}$ on $K(j)$, $j\in \Bbb{N}$ satisfies all the axioms for an
operad except that the substitution (or composition) map $\gamma_{K}$
are only partially defined. Thus $K$ is an partial operad.

For any two element of the moduli space of spheres with
tubes, in general
the $i$-th puncture of the first element might not be able to
be sewn with the $0$-th puncture of the second element.
But {}from the definition of sewing operation, we see that after rescaling the
$i$-th puncture of the first element or the $0$-th puncture
of the second  element,  the
$i$-th puncture of the first element can always
be sewn with the $0$-th puncture of the second element. So we see that
though  $K$ is partial, it is {\it rescalable}, that is, after rescaling,
the substitution map is always defined. Since all rescalings of a local
coordinate form a group isomorphic to $\Bbb{C}^{\times}$, $K$ is an example
of {\it $\Bbb{C}^{\times}$-rescalable partial operad} and $\Bbb{C}^{\times}$ is
 the {\it rescaling group} of $K$.

For any $n\in \Bbb{N}$, $K(n)$ is an infinite-dimensional complex
manifold. To be more precise, $K(n)$ is a complex (LB)-manifold, that is,
a manifold modeled on an (LB)-space over $\Bbb{C}$
(the strict inductive limit of an increasing
sequence of Banach spaces over $\Bbb{C}$) such that the transition maps
are complex analytic (see, for example, \cite{K} for the definition of
(LB)-space and  Appendix C of \cite{H4} for the precise definition of
complex (LB)-manifold).
It is proved in \cite{H1} and \cite{H4} that the sewing operation and
consequently the substitution map are
analytic and even algebraic in a certain sense.
So $K$ is an analytic $\Bbb{C}^{\times}$-rescalable partial operad.

Note that in the operadic formulation of the notion of associative algebra,
the nontrivial element $\sigma_{12}$ of $S_{2}$ generate the
operad $\{S_{j}\}_{j\in \Bbb{N}}$ for associative algebras. In fact, the
associativity of the product in an associative algebra is the reflection
of a property of $\sigma_{12}$. This element $\sigma_{12}$ is  called
an {\it associative element} of the operad $\{S_{j}\}_{j\in \Bbb{N}}$ and
$\{S_{j}\}_{j\in \Bbb{N}}$ is an example of {\it associative operads}.

For any $z\in \Bbb{C}^{\times}$, let $P(z)\in K(2)$ be
the conformal equivalence class of the sphere
$\Bbb{C}\cup\{\infty\}$ with $\infty$ the negatively oriented puncture, $z$
and $0$ the first and second positively oriented punctures,
respectively, and with the standard local coordinates vanishing at
these punctures.
Then $P(z)$ together with $K(0)$ and $K(1)$ also generates
the partial operad $K$ and $P(z)$ also has a property similar to
that of $\sigma_{12}$ (see \cite{H4}).
We call $P(z)$ an {\it associative element} of
 $K$ and $K$ is an example of
{\it associative analytic $\Bbb{C}^{\times}$-rescalable partial operads}.

In general, vertex operator algebras have nonzero central charges.
To incorporate
central charges, we need determinant line bundles and their complex powers.
For any $n\in \Bbb{N}$, there is a holomorphic line bundle $\tilde{K}(n)$
over $K(n)$ called the {\it determinant line bundle}.
The family $\tilde{K}=\{\tilde{K}(n)\}_{n\in \Bbb{N}}$ is an
associative analytic $\Bbb{C}^{\times}$-rescalable partial operad such
that the projections {}from $\tilde{K}(n)$ to $K(n)$ for all $n\in \Bbb{N}$
give
a morphism of associative analytic $\Bbb{C}^{\times}$-rescalable
partial operads and such that when restricted to the fibers,
the substitution maps are linear maps. For any $c\in \Bbb{C}$ and
$n\in \Bbb{N}$,
the complex power $\tilde{K}^{c}=\{\tilde{K}^{c}(n)\}_{n\in \Bbb{N}}$
of $\tilde{K}$ is a well-defined associative
analytic $\Bbb{C}^{\times}$-rescalable partial operad such that
for any $n\in \Bbb{N}$, $\tilde{K}^{c}(n)$ is a holomorphic line
 bundle over $K(n)$ equal to $\tilde{K}(n)$, such
that the projections {}from $\tilde{K}^{c}(n)$ to
$K(n)$ for all $n\in \Bbb{N}$ give
a morphism of associative analytic $\Bbb{C}^{\times}$-rescalable
partial operads and such that when restricted to the fibers,
the substitution maps are linear maps.
For detailed descriptions of determinant line bundles and their complex
powers, see \cite{H4}.

To define algebras over a $\Bbb{C}^{\times}$-rescalable partial operad,
we need to define a $\Bbb{C}^{\times}$-rescalable partial operad
constructed {}from a vector space. Let $V=\coprod_{n\in \Bbb{R}} V_{(n)}$ be a
$\Bbb{R}$-graded vector space such that $\dim V_{(n)}<\infty$ for
all $n\in \Bbb{R}$ and $V'=\coprod_{n\in \Bbb{R}}V_{(n)}$.
For any $n\in \Bbb{N}$,
let $\cal{H}_{V}(n)$ be the space of linear maps {}from
$V^{\otimes n}$ to $\overline{V}=\prod_{n\in \Bbb{R}}V_{(n)}=V^{\prime *}$
where $^{\prime}$ and $^{*}$ denote the functors of taking
restricted duals and duals of $\Bbb{R}$-graded vector spaces.
Since images of elements of $\cal{H}_{V}(n)$ is in general in $\overline{V}$
but not in $V$, the substitutions
are not defined in general. For example, for $f\in \cal{H}_{V}(2)$,
$g_{1}, g_{2}\in \cal{H}_{V}(1)$ and $v_{1}, v_{2}\in V$,
$f(g_{1}(v_{1}), g_{2}(v_{2}))$ is not defined in general. For any $n\in
\Bbb{C}$, let $P_{n}$ be the projection {}from $V$ to $V_{(n)}$.
Then for any
$m, n\in \Bbb{C}$, $f(P_{m}(g_{1}(v_{1})), P_{n}(g_{2}(v_{2})))$ is an element
of $\overline{V}=V^{\prime *}$. If
for any $v_{1}, v_{2}\in V$ and $v'\in V'$, the series
$\langle v', \sum_{m, n\in \Bbb{C}}f(P_{m}(g_{1}(v_{1})),
P_{n}(g_{2}(v_{2})))$
is absolutely convergent, then we obtain an element of $\cal{H}_{V}(2)$.
In this way, we obtain a partial map {}from $\cal{H}_{V}(2)\times
\cal{H}_{V}(1)\times \cal{H}_{V}(1)$ to $\cal{H}_{V}(2)$.
In general we can define partial substitution maps for the family
$\{\cal{H}_{V}(n)\}_{n\in \Bbb{N}}$ in this way.
There is also an identity $I_{V}\in \cal{H}_{V}(1)$ which is
the identity map {}from $V$ to $V\subset \overline{V}$. The symmetric group
$S_{n}$ obviously acts on $\cal{H}_{V}(n)$.
In general, however,  the substitution maps
do not satisfy the associativity even when both sides exist. We call
a structure like $\cal{H}_{V}$
a {\it partial pseudo-operad}. In particular, partial operads are
partial pseudo-operads. {\it Morphisms} between partial pseudo-operads are
defined in the obvious way.

For any $c\in \Bbb{C}$, a pseudo-algebra over $\tilde{K}^{c}$ is a
$\Bbb{Z}$-graded vector space $V$ and a morphism of partial pseudo-operads
{}from $\tilde{K}^{c}$ to $\cal{H}_{V}$. An algebra over $\tilde{K}^{c}$
is a pseudo-algebra $V$ over $\tilde{K}^{c}$ such that the image of the
morphism {}from $\tilde{K}^{c}$ to $\cal{H}_{V}$ is a partial operad.

A vertex associative algebra of central charge $c$ is a pseudo-algebra
$V$
over $\tilde{K}^{c/2}$ such that $V_{(n)}=0$ for $n$ sufficiently small and
the morphism {}from $\tilde{K}^{c/2}$ to $\cal{H}_{V}$ is meromorphic in a
certain sense (see \cite{HL1}, \cite{HL2} and \cite{H4} for the precise
definition). It turns out that the condition that
$V_{(n)}=0$ for $n$ sufficiently small and the meromorphicity of the
pseudo-algebra imply that $V$ is actually an algebra over $\tilde{K}^{c/2}$.

The following theorem announced in \cite{HL1}--\cite{HL2} and proved
in \cite{H4} is the main result of the operadic formulation of the
notion of vertex operator algebra:

\begin{theo}
The category of vertex operator algebra of central charge $c$ is isomorphic
to the category of vertex associative algebra of central charge $c$.
\end{theo}

\section{intertwining operator algebras}

In this section, the properties of the algebra of intertwining
operators for a vertex operator algebra satisfying certain conditions
are summarized to formulate the notion of intertwining operator
algebra. The axioms in the definition can be relaxed or modified to
give many variants of this notion. Since we are interested in
constructing genus-zero modular functors and weakly holomorphic
conformal field theories {}from intertwining operator algebras, we shall
only discuss the version introduced below in this paper.

Let $A$ be an $n$-dimensional commutative associative algebra over
$\Bbb{C}$. Then
for any basis $\cal{A}$ of $A$, there are structure constants $\cal{
N}_{a_{1}a_{2}}^{a_{3}}\in \Bbb{C}$,
$a_{1}, a_{2}, a_{3}\in \cal{A}$, such that
$$a_{1}a_{2}=\sum_{a_{3}\in \cal{A}}\cal{
N}_{a_{1}a_{2}}^{a_{3}}a_{3}$$ for any $a_{1}, a_{2}\in \cal{A}$.
Assume that $A$ has a basis $\cal{A}\subset A$ containing the identity
$e\in A$ such that all the structure constants $\cal{
N}_{a_{1}a_{2}}^{a_{3}}$, $a_{1}, a_{2}, a_{3}\in \cal{A}$, are in
$\Bbb{N}$.
Note that in this case, for any $a_{1}, a_{2}\in \cal{A}$,
$$\cal{N}_{ea_{1}}^{a_{2}}=\delta_{a_{1},
a_{2}}=\left\{\begin{array}{ll}
1&a_{1}=a_{2},\\0& a_{1}\ne a_{2}.\end{array}\right.$$
The commutativity and the associativity of $A$ give the following
identities:
\begin{eqnarray*}
\cal{N}_{a_{1}a_{2}}^{a_{3}}&=&\cal{N}_{a_{2}a_{1}}^{a_{3}}\\
\sum_{a\in \cal{A}}\cal{N}_{a_{1}a_{2}}^{a}\cal{
N}_{aa_{3}}^{a_{4}}&=&
\sum_{a\in \cal{A}}\cal{N}_{a_{1}a}^{a^{4}}\cal{
N}_{a_{2}a_{3}}^{a},
\end{eqnarray*}
for $a_{1}, a_{2}, a_{3}, a_{4}\in \cal{A}$.

For a vector space
$W=\coprod_{a\in \cal{A}, n\in \Bbb{R}}W^{a}_{(n)}$ doubly graded
by $\Bbb{R}$ and $\cal{A}$, let
\begin{eqnarray*}
W_{(n)}&=&\coprod_{a\in \cal{A}}W^{a}_{(n)}\\
W^{a}&=&\coprod_{n\in \Bbb{R}}W^{a}_{(n)}.
\end{eqnarray*}
Then
$$W=\coprod_{n\in \Bbb{R}}W_{(n)}=\coprod_{a\in \cal{A}}W^{a}.$$
The {\it graded dual of $W$} is the vector space
$W'=\coprod_{a\in \cal{A}, n\in \Bbb{R}}(W^{a}_{(n)})^{*}$ doubly graded
by $\Bbb{R}$ and $\cal{A}$, where for any vector space $V$, $V^{*}$
denotes the dual space of $V$. We shall denote the canonical pairing
between $W'$ and $W$ by $\langle\cdot, \cdot\rangle_{W}$.

For a vector space $W$ and a formal variable $x$, we shall in this
paper denote the space of all formal sums of the form $\sum_{n\in
\Bbb{R}}w_{n}x^{n}$ by $W\{x\}$. Note that we only allow real powers
of $x$.  For any $z\in \Bbb{C}$, we shall always choose $\log z$ so
that
$$
\log z=\log |z|+i\arg z\;\;\mbox{\rm with}\;\;0\le\arg z<2\pi.
$$

\begin{defi}
{\rm  An {\it intertwining operator algebra of central charge
$c\in \Bbb{C}$} consists of the following data:

\begin{enumerate}

\item A finite-dimensional commutative associative algebra $A$ and a
basis $\cal{A}$ of $A$ containing the identity $e\in A$ such that all
the structure constants $\cal{ N}_{a_{1}a_{2}}^{a_{3}}$, $a_{1},
a_{2}, a_{3}\in \cal{A}$, are in $\Bbb{N}$.

\item  A vector space
$$W=\coprod_{a\in
\cal{A}, n\in \Bbb{R}}W^{a}_{(n)},\; \text{for}\; w\in W^{a}_{(n)},\;
n=\wt w,\; a=\clr w$$
doubly graded
by $\Bbb{R}$ and $\cal{A}$
 (graded by {\it weight} and by {\it color}, respectively).

\item For each triple $(a_{1}, a_{2}, a_{3})\in \cal{A}\times
\cal{A}\times \cal{A}$, an $\cal{N}_{a_{1}a_{2}}^{a_{3}}$-dimensional subspace
$\cal{V}_{a_{1}a_{2}}^{a_{3}}$ of the vector space of all linear maps
$W_{a_{1}}\otimes W_{a_{2}}\to W_{a_{3}}\{x\}$, or equivalently,
 an $\cal{N}_{a_{1}a_{2}}^{a_{3}}$-dimensional vector space
$\cal{V}_{a_{1}a_{2}}^{a_{3}}$ whose elements are linear maps
\begin{eqnarray*}
\cal{Y}: W_{a_{1}}&\to& \mbox{Hom}(W_{a_{2}}, W_{a_{3}})\{x\}\\
w&\mapsto &\cal{Y}(w, x)
=\sum_{n\in \Bbb{R}}\cal{Y}_{n}(w)x^{-n-1} \;(\mbox{where}\;
\cal{Y}_{n}(w)\in \mbox{End}\;W).
\end{eqnarray*}

\item Two distinguished vectors $\bold{1}\in W^{e}$ ({\it the vacuum})
and $\omega\in W^{e}$ ({\it the Virasoro element}).

\end{enumerate}

\noindent These data satisfy the
following conditions for $a_{1}, a_{2}, a_{3}, a_{4}, a_{5},
a_{6}\in \cal{A}$,
$w_{i}\in W_{a_{i}}$, $i=1, 2, 3$, and $w'\in W'_{a_{4}}$:
\begin{enumerate}

\item The {\it grading-restriction conditions}:
\begin{eqnarray*}
&\dim W^{a}_{(n)}< \infty\;\mbox{for}\; n\in \Bbb{Z}, a\in \cal{
A},&\\
&W^{a}_{(n)}=0\;\mbox{for $n$
sufficiently small and for all}\; a\in \cal{A}\;,
\end{eqnarray*}

\item The {\it single-valuedness condition}: for any $\cal{Y}\in \cal{
V}_{ea_{1}}^{a_{1}}$,
$$\cal{Y}(w_{1}, x)\in \mbox{Hom}(W_{a_{1}},
W_{a_{1}})[[x, x^{-1}]].$$

\item The {\it lower-truncation property for vertex operators}: for any
$\cal{Y}\in \cal{V}_{a_{1}a_{2}}^{a_{3}}$,
$\cal{Y}_{n}(w_{1})w_{2}=0$ for $n$ sufficiently large.

\item The {\it identity property}: for any $\cal{Y}\in \cal{
V}_{ea_{1}}^{a_{1}}$, there is $\lambda_{\cal{Y}}\in \Bbb{C}$ such that
$\cal{Y}(\bold{1}, x)=\lambda_{\cal{Y}}I_{W_{a_{1}}}$ where
$I_{W_{a_{1}}}$ on the right is the identity
operator on $W_{a_{1}}$.

\item The {\it creation property}: for any $\cal{Y}\in \cal{
V}_{a_{1}e}^{a_{1}}$, there is $\mu_{\cal{Y}}\in \Bbb{C}$ such that $\cal{
Y}(w_{1}, x)\bold{1}\in W[[x]]$ and $\lim_{x\to 0}\cal{Y}(w_{1},
x)\bold{1}=\mu_{\cal{Y}}w_{1}$ (that is, $\cal{Y}(w_{1},
x)\bold{1}$ involves only nonnegative
integral powers of $x$ and the constant term is $\mu_{\cal{Y}}w_{1}$).

\item The {\it convergence properties}: for any $m\in \Bbb{Z}_{+}$,
$a_{i}, b_{i}, c_{i}\in \cal{A}$, $w_{i}\in W^{a_{i}}$, $\cal{Y}_{i}\in \cal{
V}_{a_{i}\;b_{i+1}}^{c_{i}}$, $i=1, \dots, m$, $w'\in (W^{c_{1}})'$ and
$w\in W^{b_{m}}$, the series
$$\langle w', \cal{Y}_{1}(w_{1}, x_{1})\cdots\cal{Y}_{m}(w_{m},
x_{m})w\rangle_{W^{c_{1}}}\mbar_{x^{n}_{i}=e^{n\log z_{i}},
i=1, \dots, m, n\in \Bbb{R}}$$
is absolutely convergent when $|z_{1}|>\cdots >|z_{m}|>0$, and for
any $\cal{Y}_{1}\in \cal{
V}_{a_{1}a_{2}}^{a_{5}}$ and $\cal{Y}_{2}\in \cal{
V}_{a_{5}a_{3}}^{a_{4}}$, the series
$$\langle w', \cal{Y}_{2}(\cal{Y}_{1}(w_{1}, x_{0})w_{2},
x_{2})w_{3}\rangle_{W}\mbar_{x^{n}_{0}=e^{n\log z_{1}-z_{2}},
x^{n}_{2}=e^{n\log z_{2}}, n\in \Bbb{R}}$$ is absolutely convergent when
$|z_{2}|>|z_{1}-z_{2}|>0$.

\item The {\it associativity}: for any $\cal{Y}_{1}\in \cal{
V}_{a_{1}a_{5}}^{a_{4}}$ and $\cal{
V}_{a_{2}a_{3}}^{a_{5}}$, there exist $\cal{Y}^{a}_{3}\in \cal{
V}_{a_{1}a_{2}}^{a}$ and $\cal{Y}^{a}_{4}\in \cal{
V}_{aa_{3}}^{a_{4}}$ for all $a\in \cal{A}$
such that the (multi-valued) analytic function
$$\langle w',
\cal{Y}_{1}(w_{1}, x_{1})\cal{Y}_{2}(w_{2},
x_{2})w_{3}\rangle_{W}\mbar_{x_{1}=z_{1},
x_{2}=z_{2}}$$
on $\{(z_{1}, z_{2})\in \Bbb{C}\times \Bbb{C}\;|\; |z_{1}|>|z_{2}|>0\}$
and the (multi-valued) analytic function
$$\sum_{a\in \cal{A}}\langle w', \cal{Y}^{a}_{4}(\cal{Y}^{a}_{3}(w_{1},
x_{0})w_{2}, x_{2})w_{3}\rangle_{W}\mbar_{x_{0}=z_{1}-z_{2},
x_{2}=z_{2}}$$
on
$\{(z_{1}, z_{2})\in \Bbb{C}\times \Bbb{C}\;|\;
|z_{2}|>|z_{1}-z_{2}|>0\}$ are equal on $\{(z_{1}, z_{2})\in
\Bbb{C}\times \Bbb{C}\;|\;|z_{1}|> |z_{2}|>|z_{1}-z_{2}|>0\}$.

\item The {\it Virasoro algebra relations}: Let $Y$ be the
element of $\cal{V}_{ea_{1}}^{a_{1}}$ such that $Y(\bold{1},
x)=I_{W_{a_{1}}}$ and let $Y(\omega, x)=\sum_{n\in
\Bbb{Z}}L(n)x^{-n-2}$. Then
$$[L(m), L(n)]=(m-n)L(m+n)+\frac{c}{12}(m^{3}-m)\delta_{m+n, 0}$$
for $m, n\in \Bbb{Z}$.

\item The {\it $L(0)$-grading property}:
$L(0)w=nw=(\wt w)w$ for $n\in \Bbb{R}$
and $w\in W_{(n)}$.

\item The {\it $L(-1)$-derivative property}: For any $\cal{Y}\in \cal{
V}_{a_{1}a_{2}}^{a_{3}}$,
$$\frac{d}{dx}\cal{Y}(w_{1}, x)=\cal{Y}(L(-1)w_{1}, x).$$

\item The {\it skew-symmetry}: There is a linear map $\Omega$ {}from
$\cal{V}_{a_{1}a_{2}}^{a_{3}}$ to $\cal{V}_{a_{2}a_{1}}^{a_{3}}$
such that for any $\cal{Y}\in \cal{V}_{a_{1}a_{2}}^{a_{3}}$,
$$\cal{Y}(w_{1}, x)w_{2}=e^{xL(-1)}(\Omega(\cal{Y}))(w_{2},
y)w_{1}\mbar_{y^{n}=e^{in\pi}x^{n}}.$$
\end{enumerate}}
\end{defi}

We shall denote the intertwining operator algebra defined above by $(W,
\cal{A},  \{\cal{
V}_{a_{1}a_{2}}^{a_{3}}\}, \bold{1}, \omega)$ or simply by $W$. The
commutative associative algebra $A$ is called the {\it Verlinde
algebra} or the {\it fusion algebra of $W$}. The linear maps in
$\cal{
V}_{a_{1}a_{2}}^{a_{3}}$ are called {\it intertwining operators of type
${a_{3}\choose a_{1}a_{2}}$}.

\begin{rema}\label{commu}
In the definition above,  the second convergence property can be
derived {}from the first one using the skew-symmetry. We include the second
convergence property in the definition in order to state the associativity
without proving this fact. The associativity or the skew-symmetry can
be replaced by the following {\it commutativity}:
for any $\cal{Y}_{1}\in \cal{
V}_{a_{1}a_{5}}^{a_{4}}$ and $\cal{Y}_{2}\in \cal{
V}_{a_{2}a_{3}}^{a_{5}}$, there exist $\cal{Y}^{a}_{5}\in \cal{
V}_{a_{2}a}^{a_{4}}$ and $\cal{Y}^{a}_{6}\in \cal{
V}_{a_{1}a_{3}}^{a}$ for all $a\in \cal{A}$
such that the (multi-valued) analytic function
$$\langle w',
\cal{Y}_{1}(w_{1}, x_{1})\cal{Y}_{2}(w_{2},
x_{2})w_{3}\rangle_{W}\lbar_{x_{1}=z_{1},
x_{2}=z_{2}}$$
on $\{(z_{1}, z_{2})\in \Bbb{C}\times \Bbb{C}\;|\; |z_{1}|>|z_{2}|>0\}$
and the (multi-valued) analytic function
$$\langle w',
\cal{Y}_{5}(w_{2}, x_{1})\cal{Y}_{6}(w_{1},
x_{2})w_{3}\rangle_{W}\lbar_{x_{1}=z_{1},
x_{2}=z_{2}}$$
on $\{(z_{1}, z_{2})\in \Bbb{C}\times \Bbb{C}\;|\; |z_{2}|>|z_{1}|>0\}$
are analytic extensions of each other.
\end{rema}

\begin{rema}\label{f-b}
The associativity gives a linear map {}from $\cal{
V}_{a_{1}a_{5}}^{a_{4}}\otimes \cal{ V}_{a_{2}a_{3}}^{a_{5}}$ to
$\oplus_{a\in \cal{A}}\cal{ V}_{a_{1}a_{2}}^{a}\otimes \cal{
V}_{aa_{3}}^{a_{4}}$ for any $a_{1}, a_{2}, a_{3}, a_{4}, a_{5}\in
\cal{A}$.  Thus we obtain a linear map {}from $\oplus_{a\in
\cal{A}}\cal{ V}_{a_{1}a}^{a_{4}}\otimes \cal{ V}_{a_{2}a_{3}}^{a}$ to
$\oplus_{a\in \cal{A}}\cal{ V}_{a_{1}a_{2}}^{a}\otimes \cal{
V}_{aa_{3}}^{a_{4}}$ for any $a_{1}, a_{2}, a_{3}, a_{4}\in \cal{A}$.
It is easy to show that these linear maps are in fact linear
isomorphisms.  These linear isomorphisms are called the {\it fusing
isomorphisms} and the associated matrices under any basis are called
{\it fusing matrices}.  Similarly,
the commutativity in Remark \ref{commu} gives a linear isomorphism
{}from $\oplus_{a\in \cal{A}}\cal{ V}_{a_{1}a}^{a_{4}}\otimes \cal{
V}_{a_{2}a_{3}}^{a}$ to $\oplus_{a\in \cal{A}}\cal{
V}_{a_{2}a}^{a_{4}}\otimes \cal{ V}_{a_{1}a_{3}}^{a}$ for any $a_{1},
a_{2}, a_{3}, a_{4}\in \cal{A}$.  These linear isomorphisms are
called the {\it braiding isomorphisms} or the associated matrices
under any basis are called
{\it braiding matrices}.
\end{rema}

In the rest of this section, we assume that the reader is familiar
with the notions of abelian intertwining algebra, module for a
vertex operator algebra, and intertwining operator among three modules
for a vertex operator algebra. We also assume that the reader knows
the basic concepts and results in the representation theory of vertex operator
algebras, for example, the rationality of vertex operator algebras and
the conditions to use the tensor product theory for modules for a
vertex operator algebra. See \cite{DL},
\cite{FLM}, \cite{FHL}, \cite{HL3}--\cite{HL6} and \cite{H3}.

It is easy to verify the following:

\begin{propo}
An abelian intertwining algebra in the sense of Dong and Lepowsky
\cite{DL} satisfying in addition the grading-restriction conditions is an
intertwining operator algebra whose Verlinde algebra
is the group algebra of the
abelian group associated with the abelian intertwining algebra.
In particular, vertex operator algebras are
intertwining operator algebras.
\end{propo}

The results in \cite{FHL}, \cite{HL3}, \cite{HL4}, and \cite{HL6}
imply the following
result:

\begin{theo}
Let $V$ be a rational vertex operator algebra satisfying the
conditions to use the tensor product theory for $V$-modules,
$\cal{A}=\{a_{1}, \dots, a_{m}\}$
the set of all equivalence classes of irreducible $V$-modules, and
$W^{a_{1}}, \dots, W^{a_{m}}$ representatives of
$a_{1},\dots, a_{m}$, respectively.
Then there is a natural commutative associative
algebra structure on the vector space $A$ spanned by $\cal{A}$ such
that $W=\coprod_{i=1}^{m}W^{a_{i}}$ has a natural structure of
intertwining operator algebra whose Verlinde algebra is $A$.
\end{theo}

We know that minimal Virasoro vertex operator
algebras are rational (see \cite{DMZ} and \cite{Wa}).
In \cite{H5}, it is proved that any vertex operator algebra containing
a vertex operator subalgebra isomorphic to a
tensor product algebra of  minimal
Virasoro vertex operator algebras satisfies
the conditions to use the tensor product theory. The
proof uses the representation theory of the Virasoro algebra and the
Belavin-Polyakov-Zamolodchikov equations (see \cite{BPZ}) for minimal
models.  We also know that the vertex operator algebras associated to
Wess-Zumino-Novikov-Witten models (WZNW models) are rational (see
\cite{FZ}).  It is easy to see that if the representation theory of
the Virasoro algebra and the Belavin-Polyakov-Zamolodchikov equations
for minimal models are replaced by the representation theory of affine
Lie algebras and the Knizhnik-Zamolodchikov equations (see \cite{KZ})
for WZNW models, respectively, the methods in \cite{H5} works for the
vertex operator algebras associated to WZNW models. Thus the vertex
operator algebras associated to WZNW models also satisfy the
conditions to use the tensor product theory.  So we have:

\begin{corol}\label{mi-wznw}
Let $V$ be a rational vertex operator algebra containing a
vertex operator subalgebra isomorphic to a tensor product algebra of  minimal
Virasoro vertex operator algebras
or a vertex operator algebras associated to a WZNW model,
$\cal{A}=\{a_{1}, \dots, a_{m}\}$
the set of all equivalence classes of irreducible $V$-modules, and
$W^{a_{1}}, \dots, W^{a_{m}}$ representatives of
$a_{1},\dots, a_{m}$, respectively.
Then there is a natural commutative associative
algebra structure on the vector space $A$ spanned by $\cal{A}$ such
that $W=\coprod_{i=1}^{m}W^{a_{i}}$ has a natural structure of
intertwining algebra whose Verlinde algebra is $A$.
\end{corol}

It is easy to see that in general the Verlinde algebras for
intertwining operator algebras obtained {}from these vertex operator
algebras are not group algebras. Thus these intertwining operator
algebras are not abelian intertwining algebras. So we do have many
examples of intertwining operator algebras which are not abelian
intertwining algebras.

\section{Genus-zero modular functors}

In the definition of intertwining operator algebra, one of the data is
the collections of the vector
spaces $\cal{V}^{a_{3}}_{a_{1}a_{2}}$, $a_{1}, a_{2}, a_{3}\in
\cal{A}$. {}From these vector spaces, without using the properties of
the operators in them, we can construct a geometric object. The
properties of these geometric objects are in fact the axioms in the
notions of genus-zero modular functors and rational genus-zero modular
functor.

We first give the definition of genus-zero modular functor.

\begin{defi}
{\rm A {\it genus-zero modular functor} is
an analytic $\Bbb{C}^{\times}$-rescalable partial
operad $\cal{M}$ together with a morphism $\pi:
\cal{M}\to K$ of $\Bbb{C}^{\times}$-rescalable partial
operads satisfying the following axioms:

\begin{enumerate}

\item For any $j\in \Bbb{N}$ the triple $(\cal{M}(j), K(j), \pi)$ is a
finite-rank holomorphic vector bundle over $K(j)$.

\item For any $\cal{Q}\in \cal{M}(k)$, $\cal{Q}_{1}\in
\cal{M}(j_{1}), \dots, \cal{Q}_{k}\in \cal{M}(j_{k})$, $k,
j_{1}, \dots, j_{k}\in \Bbb{N}$, the substitution
$\gamma_{\cal{M}}(\cal{Q}; \cal{Q}_{1},\dots, \cal{Q}_{k})$ in
$\cal{M}$ exists if (and only if)
$$\gamma_{K}(\pi(\cal{Q});
\pi(\cal{Q}_{1}), \dots, \pi(\cal{Q}_{k}))$$
exists;

\item Let $Q\in K(k)$, $Q_{1}\in K(j_{1}), \dots, Q_{k}\in K(j_{k})$,
$k, j_{1}, \dots, j_{k}\in \Bbb{N}$, such that $\gamma(Q; Q_{1},
\dots, Q_{k})$ exists.  The map {}from the Cartesian product of the
fibers over $Q$, $Q_{1}, \dots, Q_{k}$ to the fiber over
$\gamma_{K}(Q; Q_{1}, \dots, Q_{k})$ induced {}from the substitution map
of $\cal{M}$ is multilinear and gives an isomorphism {}from the tensor
product of the fibers over $Q$, $Q_{1}, \dots, Q_{k}$ to the fiber
over $\gamma_{K}(Q; Q_{1}, \dots, Q_{k})$.

\item The actions of the symmetric groups on $\cal{M}$ are isomorphisms
of holomorphic vector bundles covering the actions of the
symmetric groups on $K$.

\end{enumerate}}
\end{defi}

{\it Homomorphisms} and {\it isomorphisms}
{}from genus-zero modular functors to genus-zero modular functors are defined
in the obvious way.

Let $\cal{M}$ be a genus-zero modular functor.  For any $k\in
\Bbb{Z}_{+}$, $j_{1},\dots, j_{k}\in \Bbb{N}$, let $\cal{M}_{i}(k)$,
$\cal{M}_{i}(j_{1}), \dots, \cal{M}_{i}(j_{k})$ be vector subbundles of
$\cal{M}(k)$, $\cal{M}(j_{1}), \dots, \cal{M}(j_{k})$, respectively,
for $i=1, 2$, and $P_{i}(k)$, $P_{i}(j_{1}), \dots, P_{i}(j_{k})$ the
projections {}from $\cal{M}(k)$, $\cal{M}(j_{1}), \dots, \cal{M}(j_{k})$
to $\cal{M}_{i}(k)$,
$\cal{M}_{i}(j_{1}), \dots, \cal{M}_{i}(j_{k})$, respectively, for
$i=1, 2$. Then $\gamma_{\cal{M}}\circ (P_{1}(k)\times
P_{1}(j_{1})\cdots \times P_{1}(j_{k}))$ and $\gamma_{\cal{M}}\circ
(P_{2}(k)\times P_{2}(j_{1})\cdots \times P_{2}(j_{k}))$ are both
homomorphisms of vector bundles
{}from
$\cal{M}(k)\times \cal{M}(j_{1})\times
\cdots\times\cal{M}(j_{k})$ to
$\cal{M}(j_{1}+\cdots +j_{k})$.
So the sum of $\gamma_{\cal{M}}\circ (P_{1}(k)\times
P_{1}(j_{1})\cdots \times P_{1}(j_{k}))$ and $\gamma_{\cal{M}}\circ
(P_{2}(k)\times P_{2}(j_{1})\cdots \times P_{2}(j_{k}))$ is
well-defined. Similarly, sums of more than two such
homomorphisms is also well-defined.

We next give the definition of  rational genus-zero modular functor:

\begin{defi}
A {\it rational genus-zero modular functor} is a genus-zero modular
functor $\cal{M}$ and a finite set $\cal{A}$ satisfying the
following axioms:

\begin{enumerate}

\setcounter{enumi}{4}

\item For any $n\in \Bbb{N}$ and $a_{0}, a_{1}, \dots, a_{n}\in \cal{
A}$, there are finite-rank holomorphic vector bundles $\cal{
M}_{a_{1}\cdots a_{n}}^{a_{0}}(n)$ over $K(n)$ such that $\cal{
M}(n)=\oplus_{a_{0}, a_{1}, \dots, a_{n}\in \cal{A}}\cal{
M}_{a_{1}\cdots a_{n}}^{a_{0}}(n)$ where $\oplus$ denotes the direct
sum operation for vector bundles.

\item For any $n\in \Bbb{N}$, $b_{0}, b_{1}, \dots, b_{n}\in \cal{A}$,
let $P^{b_{0}}(n)$ be the projection {}from $\cal{M}(n)$ to
$\oplus_{a_{1}, \dots, a_{n}\in \cal{A}}\cal{ M}_{a_{1}\cdots
a_{n}}^{b_{0}}(n)$ and $P_{b_{1}\cdots b_{n}}(n)$ the projection {}from
$\cal{M}(n)$ to $\oplus_{a_{0}\in \cal{A}} \cal{ M}_{b_{1}\cdots
b_{n}}^{a_{0}}(n)$. Then for any $k\in \Bbb{Z}_{+}$, $j_{1}, \dots,
j_{k}\in \Bbb{N}$,
$$\gamma_{\cal{M}}
= \sum_{b_{1}, \dots, b_{k}\in {\cal
A}}\gamma_{\cal{M}}\circ (P_{b_{1}\cdots b_{k}}(k)
\times P^{b_{1}}(j_{1})\times \cdots P^{b_{k}}(j_{k})).$$
\end{enumerate}
\end{defi}

We shall denote the rational genus-zero modular functor just defined by
$(\cal{M}, \cal{A})$.

The simplest examples of rational genus-zero modular functors
 $\tilde{K}^{c}$, $c\in \Bbb{C}$.

Using the methods developed in \cite{H4}, we have:

\begin{theo}\label{ioa-mf}
Let $(W,
\cal{A},  \{\cal{
V}_{a_{1}a_{2}}^{a_{3}}\}, \bold{1}, \omega)$ be an intertwining operator
algebra of central charge $c$.
Then there is  a canonical rational genus-zero
modular functor $(\cal{M}_{W}, \cal{A})$ such that for any $a_{1},
a_{2}, a_{3}\in \cal{A}$,
$\cal{M}^{a_{1}}(0)=\tilde{K}^{c/2}(0)$ and the fiber of
$\cal{M}_{a_{1}a_{2}}^{a_{3}}(2)$
is isomorphic to $\cal{V}_{a_{1}a_{2}}^{a_{3}}$.
\end{theo}

In particular, by Theorem \ref{ioa-mf} and Corollary \ref{mi-wznw},
we have:

\begin{corol}
Let $V$ be a rational vertex operator algebra containing a
vertex operator subalgebra isomorphic to a tensor product algebra of  minimal
Virasoro vertex operator algebras
or a vertex operator algebras associated to a WZNW model,
$\cal{A}=\{a_{1}, \dots, a_{m}\}$
the set of all equivalence classes of irreducible $V$-modules, and
$W^{a_{1}}, \dots, W^{a_{m}}$ representatives of
$a_{1},\dots, a_{m}$, respectively.
Then there is  a canonical rational genus-zero
modular functor $(\cal{M}_{V}, \cal{A})$ such that for any $a_{1},
a_{2}, a_{3}\in \cal{A}$,
$\cal{M}^{a_{1}}(0)=\tilde{K}^{c/2}(0)$ and the fiber of
$\cal{M}_{a_{1}a_{2}}^{a_{3}}(2)$
is isomorphic to $\cal{V}_{a_{1}a_{2}}^{a_{3}}$,
where $c\in \Bbb{C}$ is the central charge of $V$ and
$\cal{V}_{a_{1}a_{2}}^{a_{3}}$ is the space of all intertwining
operators
for $V$
of type ${W^{a_{3}}\choose W^{a_{1}}W^{a_{2}}}$.
\end{corol}

Genus-zero modular functors contain almost all of the topological information
one can obtain {}from an intertwining operator algebra. To see this, we
need a result of Segal \cite{S}:

\begin{propo}
Let $\cal{M}$ be a genus-zero modular functor. Then there are canonical
projectively flat connections on the vector bundles $\cal{M}(n)$,
$n\in \Bbb{N}$, such that they are compatible with the substitution
maps for the partial operad $\cal{M}$.
\end{propo}

For any $n\in \Bbb{Z}_{+}$, we can embed the configuration space
$$F(n)=\{(z_{1}, \dots, z_{n})\in \Bbb{C}^{n}\;|\;z_{i}\ne z_{j}\;\mbox{for}\;
i\ne j\}$$ into $K(n)$ in the obvious way. Thus given any
genus-zero modular functor $\cal{M}$,  for any
$n\in \Bbb{Z}_{+}$,
$\cal{M}(n)$ is pulled back to a vector bundle over $F(n)$ and the
canonical connection on $\cal{M}(n)$ is pulled back to a connection
on the pull-back vector bundle over $F(n)$. Since the actions of symmetric
groups on $\cal{M}$
are isomorphisms of holomorphic vector bundles covering the actions of the
symmetric actions on $K$, the pull-back vector bundle and the pull-back
connection over $F(n)$ induce a vector bundle over $F(n)/S_{n}$ and
a connection over this vector bundle for any $n\in \Bbb{Z}_{+}$.
It is not difficult to prove the following:

\begin{propo}
Let $\cal{M}$ be a genus-zero modular functor.
For any $n\in \Bbb{Z}_{+}$, the connections on the pull-back vector bundle
over $F(n)$ and on the induced vector bundle over
$F(n)/S_{n}$ are flat.
\end{propo}

We know that a flat connection on a vector bundle over a connected
manifold gives a structure of a representation of the fundamental
group of the manifold to the fiber of the vector bundle at any point
on the manifold.  We also know that the braid group $B_{n}$ of $n$
strings is by definition the fundamental group of $F(n)/S_{n}$. So we
obtain:

\begin{theo}
Let $\cal{M}$ be a genus-zero modular functor, $n\in \Bbb{Z}_{+}$ and
$\cal{V}$ the fiber at any point in $F(n)\subset K(n)$ of the vector
bundle $\cal{M}(n)$. Then $\cal{V}$ has a natural structure of a
representation of the braid group $B_{n}$ of $n$ strings.
\end{theo}

In particular, {}from a vertex operator algebra containing a subalgebra
isomorphic to  a tensor product algebra of  minimal
Virasoro vertex operator algebras or {}from a
vertex operator algebra associated to WZNW models, we obtain
representations of braid groups. In the case that the vertex
operator algebras are associated to WZNW models based on the Lie group
$SU(2)$, the corresponding representations of the braid groups are the
same as those obtained by Tsuchiya and Kanie \cite{TK} and the
corresponding knot invariants are  the Jones polynomials for knots
\cite{J}.
(In fact, {}from a genus-zero modular functor, we obtain not only
representations of the braid groups, but also representations of braid
groups with twists.)

Let $\cal{M}=\{\cal{M}(n)\}_{n\in \Bbb{N}}$ be a genus-zero modular functor
and $\overline{\cal{M}}(n)$, $n\in \Bbb{N}$, the complex conjugates of the
holomorphic vector bundles $\cal{M}(n)$. We define the {\it complex conjugate
partial operad
of $\cal{M}$}
to be the sequence $\overline{\cal{M}}
=\{\overline{\cal{M}}(n)\}_{n\in \Bbb{N}}$ with the obvious partial
operad structure induced {}from that of $\cal{M}$. Let
$$\cal{M}\otimes \overline{\cal{M}}=
\{\cal{M}(n)\otimes \overline{\cal{M}}(n)\}_{n\in \Bbb{N}}$$
 where on the right-hand side, $\otimes$
is the tensor product operation for vector bundles over $K(n)$, $n\in \Bbb{N}$.
Then the partial operad structures on $\cal{M}$ and $\overline{\cal{M}}$
induce a partial operad structure on $\cal{M}\otimes \overline{\cal{M}}$.
In particular, for any $c\in \Bbb{C}$, since $\tilde{K}^{c/2}$ is
a genus-zero modular functor, we obtain partial operads
$\overline{\tilde{K}}^{c/2}$ and $\tilde{K}^{c/2}\otimes
\overline{\tilde{K}}^{c/2}$.
Using these partial operads, we have the following notion:

\begin{defi}
A genus-zero modular functor $\cal{M}$ is {\it unitary} if there is
a complex number $c\in \Bbb{C}$ and
a morphism of analytic partial operads {}from $\cal{M}\otimes
\overline{\cal{M}}$ to
$\tilde{K}^{c/2}\otimes
\overline{\tilde{K}}^{c/2}$ satisfying the following conditions:

\begin{enumerate}

\item  For any $n\in \Bbb{N}$ and $Q\in K(n)$,
the morphism {}from $\cal{M}\otimes
\overline{\cal{M}}$ to
$\tilde{K}^{c/2}\otimes
\overline{\tilde{K}}^{c/2}$ maps the fiber of $\cal{M}(n)\otimes
\overline{\cal{M}}(n)$ at $Q$ linearly
to the fiber of $\tilde{K}^{c/2}(n)\otimes
\overline{\tilde{K}}^{c/2}(n)$ at $Q$.

\item The linear map above {}from the fiber of $\cal{M}(n)\otimes
\overline{\cal{M}}(n)$ at $Q$
to the fiber of $\tilde{K}^{c/2}(n)\otimes
\overline{\tilde{K}}^{c/2}(n)=\Bbb{C}$ at $Q$ when viewed as a
linear map {}from the fiber of $\cal{M}(n)\otimes
\overline{\cal{M}}(n)$ at $Q$
to $\Bbb{C}$ gives a positive-definite Hermitian form on
the fiber of $\cal{M}(n)$ at $Q$.

\end{enumerate}
\end{defi}

The complex number $c$ is called the {\it central charge} of the
unitary
genus-zero modular functor $\cal{M}$.

Let $(W, \cal{A},  \{\cal{V}_{a_{1}a_{2}}^{a_{3}}\},
\bold{1}, \omega)$ be an intertwining
algebra of central charge $c$. Assume that there are
positive-definite Hermitian forms
on $\cal{V}_{a_{1}a_{2}}^{a_{3}}$ for all $a_{1}, a_{2}, a_{3}\in \cal{A}$.
These positive-definite Hermitian forms induce positive-definite
Hermitian forms on $\oplus_{a\in \cal{A}}\cal{
V}_{a_{1}a}^{a_{4}}\otimes \cal{
V}_{a_{2}a_{3}}^{a}$ and $\oplus_{a\in \cal{A}}\cal{
V}_{a_{1}a_{2}}^{a}\otimes \cal{
V}_{aa_{3}}^{a_{4}}$ for any $a_{1}, a_{2}, a_{3}, a_{4}\in \cal{A}$.

We have:

\begin{propo}
Let $(W, \cal{A},  \{\cal{V}_{a_{1}a_{2}}^{a_{3}}\},
\bold{1}, \omega)$ be an intertwining operator
algebra of central charge $c$. If there are  positive-definite Hermitian forms
on $\cal{V}_{a_{1}a_{2}}^{a_{3}}$ for all $a_{1}, a_{2}, a_{3}\in \cal{A}$
such that the fusing isomorphism {}from $\oplus_{a\in \cal{A}}\cal{
V}_{a_{1}a}^{a_{4}}\otimes \cal{
V}_{a_{2}a_{3}}^{a}$ to $\oplus_{a\in \cal{A}}\cal{
V}_{a_{1}a_{2}}^{a}\otimes \cal{
V}_{aa_{3}}^{a_{4}}$ pulls the induced positive-definite Hermitian form
on $\oplus_{a\in \cal{A}}\cal{
V}_{a_{1}a_{2}}^{a}\otimes \cal{
V}_{aa_{3}}^{a_{4}}$ back to that on $\oplus_{a\in \cal{A}}\cal{
V}_{a_{1}a}^{a_{4}}\otimes \cal{
V}_{a_{2}a_{3}}^{a}$ for any $a_{1}, a_{2}, a_{3}, a_{4}\in \cal{A}$,
then these positive-definite Hermitian forms give a unitary structure
to the rational modular functor $(\cal{M}_{W}, \cal{A})$.
\end{propo}

In particular, we have:

\begin{propo}
The genus-zero modular functors obtained {}from minimal Virasoro vertex
operator
algebras and {}from vertex operator algebras associated to WZNW models
are unitary.
\end{propo}

\section{Genus-zero conformal field theories}

Modular functors only reflects the geometric objects constructed {}from
the vector spaces $\cal{V}^{a_{3}}_{a_{1}a_{2}}$, $a_{1}, a_{2}, a_{3}\in
\cal{A}$, without using the properties of the operators in these spaces.
Also an intertwining operator algebra has the vacuum and the Virasoro element
which are not reflected in the definition of (rational)
genus-zero modular functor. These data and the properties they satisfied are
reflected in the following definition of genus-zero weakly holomorphic
conformal field theory:

\begin{defi}
Let $\cal{M}$ be a genus-zero modular functor. An {\it genus-zero
weakly holomorphic conformal field theory over $\cal{M}$} is
a $\Bbb{R}$-graded vector space $W$ and a morphism of partial pseudo-operads
{}from $\Bbb{M}$ to the endomorphism partial pseudo-operad $\cal{H}_{W}$
(see Section 2) satisfying the following additional axioms:

\begin{enumerate}

\item The image of the morphism {}from $\Bbb{M}$ to $\cal{H}_{W}$  is a
partial operad.

\item  The morphism {}from $\cal{M}$ to the
$\Bbb{C}^{\times}$-resclalable endomorphism partial pseudo-operad
$\cal{H}_{W}$ is
linear on the fibers of the bundles $\cal{M}(n)$, $n\in \Bbb{N}$.

\item This morphism is holomorphic.

\end{enumerate}

\noindent A genus-zero
weakly holomorphic conformal field theory over a rational genus-zero modular
functor is called a {\it rational genus-zero
holomorphic weakly conformal field theory}.
\end{defi}

Using the method developed in \cite{H4}, we have:

\begin{theo}\label{ioa-wcft}
Let $(W,
\cal{A},  \{\cal{
V}_{a_{1}a_{2}}^{a_{3}}\}, \bold{1}, \omega)$ be an intertwining operator
algebra and $(\cal{M}_{W}, \cal{A})$ the corresponding rational
genus-zero modular
functor. Then $W$ has a canonical structure of a
rational holomorphic genus-zero weakly
conformal field theory.
\end{theo}

In particular, by Theorem \ref{ioa-wcft} and Corollary \ref{mi-wznw},
we have:

\begin{corol}
Let $V$ be a rational vertex operator algebra containing a
vertex operator subalgebra isomorphic to a tensor product algebra of  minimal
Virasoro vertex operator algebras
or a vertex operator algebras associated to a WZNW model,
$\cal{A}=\{a_{1}, \dots, a_{m}\}$
the set of all equivalence classes of irreducible $V$-modules, and
$W^{a_{1}}, \dots, W^{a_{m}}$ representatives of
$a_{1}, \dots, a_{m}$, respectively.
Then there is a canonical rational genus-zero weakly
holomorphic conformal field theory structure on $W=\coprod_{i=1}^{m}W^{a_{i}}$.
\end{corol}

Let $\cal{H}_{1}$ be the Banach space of all analytic functions
on the closed
unit disk $\{z\in \Bbb{C}\;|\;|z|\le 1\}$ and $K_{\cal{H}_{1}}(n)$ for
any $n\in \Bbb{N}$ the subset of $K(n)$ consisting of conformal equivalence
classes containing spheres with $n+1$ tubes such that the inverses
of the  local coordinates
at punctures can be extended to the closed unit disk and
the images of the closed unit disks under these inverses are disjoint.
Let $K_{\cal{H}_{1}}
=\{K_{\cal{H}_{1}}(n)\}_{n\in \Bbb{N}}$. Then the sewing operation is always
defined in $K_{\cal{H}_{1}}$. Also the identity of $K$ is in $K_{\cal{H}_{1}}$
and the symmetric groups act on $K_{\cal{H}_{1}}$. So
$K_{\cal{H}_{1}}$ is a suboperad of $K$ (not just a partial suboperad of $K$).
For any $n\in \Bbb{N}$ and
$c\in \Bbb{C}$, let $\tilde{K}_{\cal{H}}^{c/2}(n)$ be the restriction of
the line bundle $\tilde{K}^{c/2}(n)$ to $K_{\cal{H}_{1}}(n)$. Then
$\tilde{K}_{\cal{H}}^{c/2}=\{\tilde{K}_{\cal{H}}^{c/2}(n)\}_{n\in \Bbb{N}}$
is a suboperad (not partial) of $\tilde{K}^{c/2}$.

Let $\overline{\tilde{K}}_{\cal{H}_{1}}^{c/2}$ be the complex conjugate
operad of
$\tilde{K}_{\cal{H}_{1}}^{c/2}$ defined in the obvious way and
$\tilde{K}_{\cal{H}_{1}}^{c/2} \otimes
\overline{\tilde{K}}_{\cal{H}_{1}}^{c/2}$
the tensor product operad of $\tilde{K}_{\cal{H}_{1}}^{c/2}$ and
$\overline{\tilde{K}}_{\cal{H}_{1}}^{c/2}$.
In general, algebras over the operad
$\tilde{K}_{\cal{H}_{1}}^{c/2} \otimes
\overline{\tilde{K}}_{\cal{H}_{1}}^{c/2}$
might not have any
topological structure. But since
$\tilde{K}_{\cal{H}_{1}}^{c/2} \otimes
\overline{\tilde{K}}_{\cal{H}_{1}}^{c/2}$
consists of infinite-dimensional
Banach manifolds, we are interested in algebras over it with topological
structures. Consider the abelian
category of Hilbert spaces over $\Bbb{C}$. Then the
tensor product operation for Hilbert spaces gives a tensor category structure
to this abelian category. Thus for any Hilbert space $H$ over $\Bbb{C}$,
we have the endomorphism operad of $H$. This endomorphism operad has
a topological structure. A {\it Hilbert algebra} over
$\tilde{K}_{\cal{H}_{1}}^{c/2} \otimes
\overline{\tilde{K}}_{\cal{H}_{1}}^{c/2}$
 is a Hilbert space $H$ over $\Bbb{C}$
and a continuous morphism of topological operads {}from
$\tilde{K}_{\cal{H}_{1}}^{c/2} \otimes
\overline{\tilde{K}}_{\cal{H}_{1}}^{c/2}$
to the endomorphism operad of $H$.

\begin{defi}
A genus-zero conformal field theory of central charge $c$
is a Hilbert algebra over
$\tilde{K}_{\cal{H}_{1}}^{c/2} \otimes
\overline{\tilde{K}}_{\cal{H}_{1}}^{c/2}$
such that the morphism {}from
$\tilde{K}_{\cal{H}_{1}}^{c/2} \otimes
\overline{\tilde{K}}_{\cal{H}_{1}}^{c/2}$
to the endomorphism operad of $H$ is linear on the fibers of
$\tilde{K}_{\cal{H}_{1}}^{c/2}(n) \otimes
\overline{\tilde{K}}_{\cal{H}_{1}}^{c/2}(n)$, $n\in \Bbb{N}$.
\end{defi}

Let $\cal{M}$ be a genus-zero modular functor and
$\cal{M}_{\cal{H}_{1}}(n)$, $n\in \Bbb{N}$, the restrictions of
$\cal{M}(n)$ to $K_{\cal{H}_{1}}(n)$. Then $\cal{M}_{\cal{H}_{1}}=
\{\cal{M}_{\cal{H}_{1}}(n)\}_{n\in \Bbb{N}}$ is a suboperad (not partial) of
$\cal{M}$. We also have the notion of {\it Hilbert algebra} over
$\cal{M}_{\cal{H}_{1}}$ defined in the obvious way.

\begin{defi}
Let $\cal{M}$ be a unitary genus-zero modular functor.
A genus-zero weakly holomorphic conformal field theory $W$ over $\cal{M}$
is {\it unitary} if there is a positive-definite
Hermitian form on $W$ such that
the morphism of partial pseudo-operads
 $\cal{M}$ to $\cal{H}_{W}$ induces
a Hilbert algebra structure on the completion $H_{W}^{h}$
($h$ means holomorphic)
of $W$ over
$\cal{M}_{\cal{H}_{1}}$.
\end{defi}

We have:

\begin{propo}
The genus-zero weakly holomorphic conformal field theories constructed
{}from the minimal Virasoro vertex operator algebras and {}from the vertex
operator algebras associated to WZNW models are unitary.
\end{propo}

Let $\cal{M}$ be a unitary genus-zero modular functor. Then for
any $n\in \Bbb{N}$ and $Q\in K(n)$, the corresponding
positive-definite Hermitian form on the fiber of $\cal{M}(n)$
at $Q$ identifies the dual space of
this fiber with the fiber of $\overline{\cal{M}}(n)$
at $Q$ and thus also identifies the dual space of
the fiber of $\overline{\cal{M}}(n)$
at $Q$ with the fiber of $\cal{M}(n)$
at $Q$. Thus the adjoint of
the morphism {}from $\cal{M}\otimes \overline{\cal{M}}$
to $\tilde{K}^{c/2}\otimes \overline{\tilde{K}}^{c/2}$ gives
a morphism {}from $(\tilde{K}^{c/2}\otimes \overline{\tilde{K}}^{c/2})^{-1}$
to $\cal{M}\otimes \overline{\cal{M}}$ where
$(\tilde{K}^{c/2}\otimes \overline{\tilde{K}}^{c/2})^{-1}
=\{(\tilde{K}^{c/2}\otimes
\overline{\tilde{K}}^{c/2})^{-1}(n)\}_{n\in \Bbb{N}}$, and
for any $n\in \Bbb{N}$, $(\tilde{K}^{c/2}\otimes
\overline{\tilde{K}}^{c/2})^{-1}(n)$ is the line bundle  whose fibers are
the duals of the fibers of $\tilde{K}^{c/2}(n)\otimes
\overline{\tilde{K}}^{c/2}(n)$. It is clear that $(\tilde{K}^{c/2}\otimes
\overline{\tilde{K}}^{c/2})^{-1}$ is canonically
isomorphic to $\tilde{K}^{c/2}\otimes
\overline{\tilde{K}}^{c/2}$. Thus we obtain a morphism {}from
$\tilde{K}^{c/2}\otimes
\overline{\tilde{K}}^{c/2}$ to $\cal{M}\otimes \overline{\cal{M}}$.

Let $W$ be a unitary genus-zero weakly holomorphic conformal field theory
over a unitary genus-zero modular functor $\cal{M}$.
Let $\overline{H}^{h}_{W}$  be the complex conjugate of
$H^{h}_{W}$ and $\cal{M}_{\cal{H}_{1}}$ the complex conjugate operad of
$\cal{M}_{\cal{H}_{1}}$.
Then $\overline{H}^{h}_{W}$ has a structure of Hilbert algebra over
$\cal{M}_{\cal{H}_{1}}$. Let $H_{W}=H^{h}_{W}\otimes \overline{H}^{h}_{W}$
 where $\otimes$
is the Hilbert space tensor product. Then
the Hilbert algebra structure on $H^{h}_{W}$ over $\cal{M}_{\cal{H}_{1}}$ and
the Hilbert algebra structure on $\overline{H}^{h}_{W}$ over
$\overline{\cal{M}}_{\cal{H}_{1}}$ induce a
Hilbert algebra structure on $H_{W}$ over $\cal{M}\otimes \overline{\cal{M}}$.
Since $\cal{M}$ is unitary, we have
a morphism {}from
$\tilde{K}^{c/2}\otimes
\overline{\tilde{K}}^{c/2}$ to $\cal{M}\otimes \overline{\cal{M}}$ by the
discussion above. Combining this morphism and the
Hilbert algebra structure  over $\cal{M}\otimes \overline{\cal{M}}$ on $H_{W}$,
we obtain:

\begin{propo}
Let $W$ be a unitary genus-zero weakly holomorphic conformal field theory
over a unitary genus-zero modular functor $\cal{M}$. Then $H_{W}$
is a genus-zero conformal field theory.
\end{propo}

In particular, we have:

\begin{corol}
Let $W$ be the unitary genus-zero weakly holomorphic conformal field theory
constructed {}from a minimal Virasoro vertex operator algebra
or {}from a vertex
operator algebras associated to WZNW model. Then $H_{W}$
is a genus-zero conformal field theory.
\end{corol}

\end{document}